# Binary Color-Coded Magic Squares: A Study of Uniqueness Under Rotation/Reflection, PCA, and LDA Analysis

Peyman Fahimi

*Dép. de chimie, Université Laval, Québec, Québec, Canada G1V 0A6.*
*Tel.: +1-(438)-229-2038, E-mail: peyman.fahimi.1@ulaval.ca.*

**Abstract**

In this paper, we study the concept of "binary color-coded magic squares" by assigning two distinct colors to the even and odd numbers within a magic square. We investigate the uniqueness of patterns within these squares using three different analytical methods, including rotation/reflection, PCA, and LDA. Our investigation covers all 880 magic squares of order 4, all 48,544 associative magic squares of order 5, and all 368,640 Franklin magic squares of order 8. Our investigation reveals striking patterns that were previously unknown in traditional magic squares, shedding light on the potential for binary color-coded magic squares to contribute to the field of mathematics.

## Introduction

Magic squares have long been a topic of fascination in mathematics, with the earliest known example dating back to ancient China [1]. A magic square is a square grid of numbers in which the sum of the numbers in each row, column, and main diagonal is the same [2–5]. While traditional magic squares follow strict rules and patterns, variations on this concept have led to the discovery of new and exciting properties within these structures [5–12].

In this paper, we focus on a modification of the traditional magic square concept by introducing a color-coding element. Specifically, we explore the concept of "binary color-coded magic squares," which assigns two distinct colors to the even and odd numbers within a magic square. This novel approach provides a new perspective on traditional magic squares and reveals previously unknown patterns and relationships.

Our investigation focuses on the uniqueness of patterns within binary color-coded magic squares, using three different analytical methods. Firstly, we analyze the squares through rotation and reflection, seeking patterns that remain unchanged under such transformations. Secondly, we use principal component analysis (PCA) [13] to identify the most significant patterns within the squares. Lastly, we apply linear discriminant analysis (LDA) [14] to classify and distinguish between different types of magic squares. PCA is a technique that identifies the most significant patterns within a dataset by transforming the data into a new set of variables that are uncorrelated with each other. LDA is a supervised learning technique that aims to classify and distinguish between different types of objects based on their attributes.

We analyze a large number of binary color-coded magic squares, including all 880 magic squares of order 4 [15,16], all 48,544 associative magic squares of order 5 [15,16], and all 368,640 Franklin magic squares of order 8 [16]. An associative magic square is a specific type of magic square that exhibits a unique property. In this type of magic square, each pair of numbers located symmetrically opposite to the center of the square adds up to the same value [5]. Franklin magic squares are also a special type of magic square with unique properties. In a Franklin magic square, each half-row or half-column of the square adds up to half of the magic constant. Additionally, the sum of each bent diagonal in the square is equal to the magic constant. Finally, the sum of each 2×2 block is equal to $4/n$ of the magic constant in which $n$ is the order of the square [16,17].

Our investigation reveals striking patterns that emerge from this novel approach, which were previously unknown in traditional magic squares. This study sheds light on the potential for color-coded magic squares to contribute to the field of mathematics and offers new insights into the nature of magic squares.

## Distinct Binary Color-Coded Magic Squares: A Rotation and Reflection Analysis

In a binary color-coded magic square, odd and even numbers are separated into two groups and assigned different colors. To identify the distinct color-coded patterns, a Python code (supplementary data) converts each of the squares into a binary matrix by assigning a value of 1 to the odd numbers and 0 to the even numbers.

Since a magic square can be rotated or reflected and still retain its symmetry and properties, the code excludes any matrices produced through these transformations. This is done to ensure that only unique patterns are identified, rather than different variations of the same pattern. Once the distinct binary matrices have been identified, the code converts them back into a pattern string to count the number of squares with each unique pattern.

Overall, we use a combination of matrix manipulation and pattern recognition techniques to identify and visualize the unique binary color-coded patterns within a specific type of magic square. Table 1 provides a summary of the unique binary matrices and the total number of distinct patterns for various types of magic squares. The table also includes the number of squares that have a particular pattern. To aid visualization, the unique patterns are illustrated visually in the table.

Table 1 shows that the 880 magic squares of order 4 can be created using only 8 unique binary magic squares. Similarly, the 48,544 associative magic squares of order 5 can be constructed with just 15 unique binary magic squares. Surprisingly, the 368,640 Franklin magic squares are created by only 6 distinct binary magic squares. Beautiful tilings can be created by plotting the binary color-coded magic squares of higher orders or repeating the patterns of lower order magic squares in a larger square grid. Figure 1 shows two examples.

*Table 1: Unique patterns of binary color-coded magic squares*

| Type of magic square | No. magic squares | No. unique patterns (UP) | Binary matrix of UP | No. squares with a UP | Unique patterns |
|---|---|---|---|---|---|
| **3×3 magic square** | 1 | 1 | 010111010 | 1 | 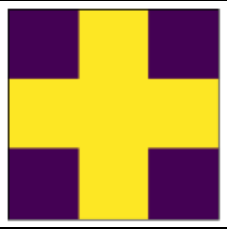 |
| **4×4 associative magic square** | 48 | 2 | 0101101010100101<br>0011110011000011 | 24<br>24 | 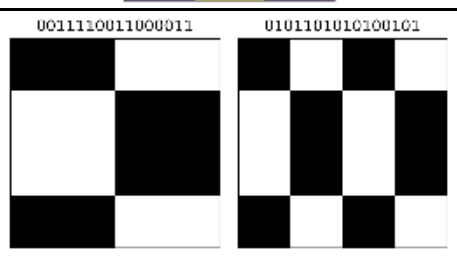 |
| **4×4 magic square** | 880 | 8 | 0011001111001100<br>0101110000111010<br>0011110011000011<br>0011110000111100<br>0011101001011100<br>0101101010100101<br>0011010110101100<br>0101101001011010 | 44<br>48<br>212<br>192<br>80<br>212<br>48<br>44 | 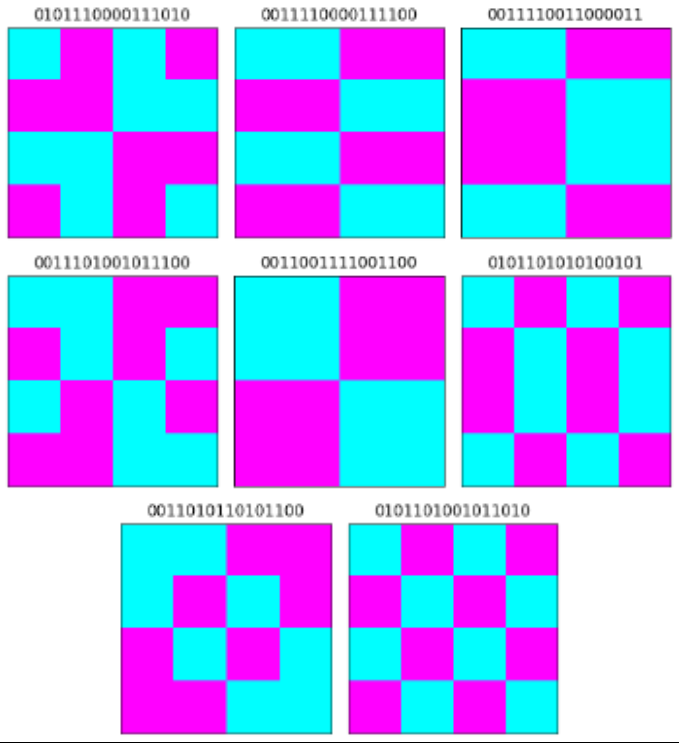 |
| **5×5 ultra magic square** | 16 | 2 | 1001101011001001101011001<br>0010010101111111010100100 | 8<br>8 | 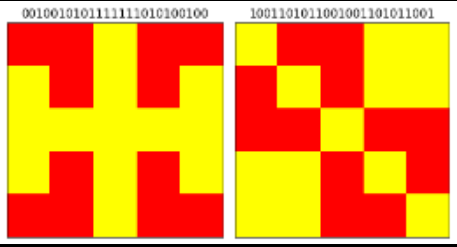 |
| **5×5 associative magic square** | 48,544 | 15 | 1010100100111110010010101<br>1001101011001001101011001<br>1001111010001000101111001<br>0001011001111111001101000<br>0001010011111111100101000<br>1001100010111110100011001<br>1001101000111110001011001<br>0010010101111111010100100<br>0000111010111110101110000<br>0000101011111111101010000<br>0101110000111110000111010<br>0011110101001001010111100<br>0101110011001001100111010<br>0101111001001001001111010<br>0010001110111110111000100 | 864<br>2180<br>2180<br>4546<br>4546<br>4546<br>4546<br>1728<br>4546<br>4546<br>4546<br>4546<br>2180<br>2180<br>864 | 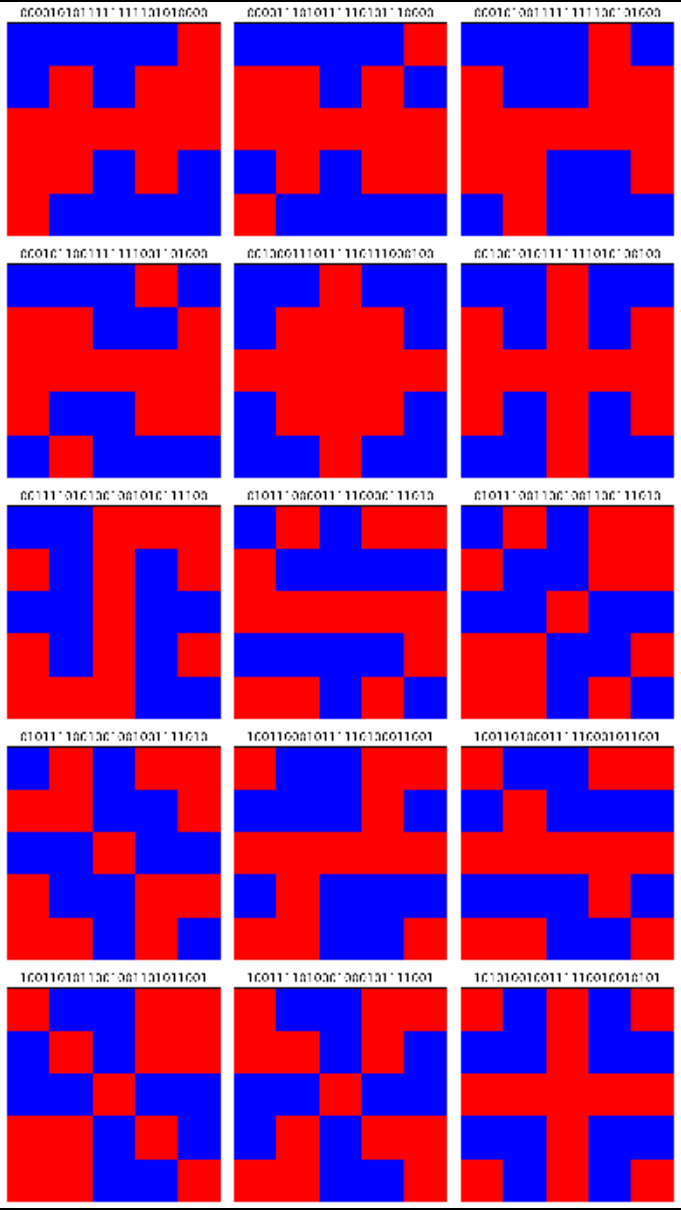 |
| **8×8 Franklin magic square** | 368,640 | 6 | 0101010110101010101010101<br>0100101010110101010101001<br>0101010101010110101010<br><br>001101101100100100110<br>110110010010011011011<br>0010010011011011001001<br><br>0101010110101010101010101<br>0100101010101010101010110<br>10101010101010100101010101<br><br>001110011100011000111<br>001110001100011100111<br>0001100011100111000110<br><br>001100111100110000110<br>011110011000011001111<br>0011000011001111001100<br><br>001111001100001100111<br>100110000110011110011<br>0000110011110011000011 | 46080<br><br>92160<br><br>46080<br><br>92160<br><br>46080<br><br>46080 | 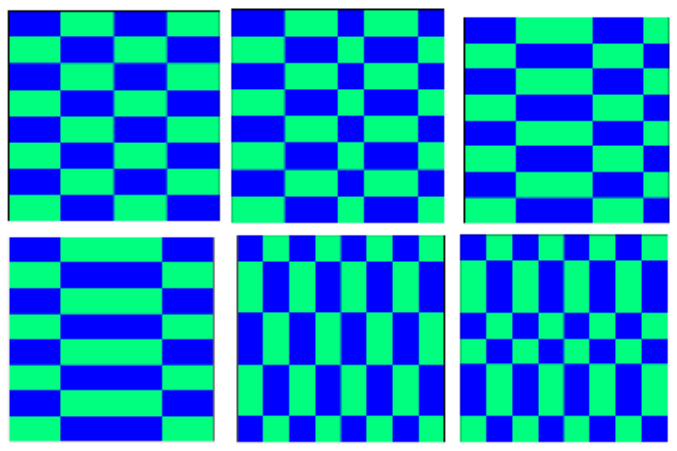 |

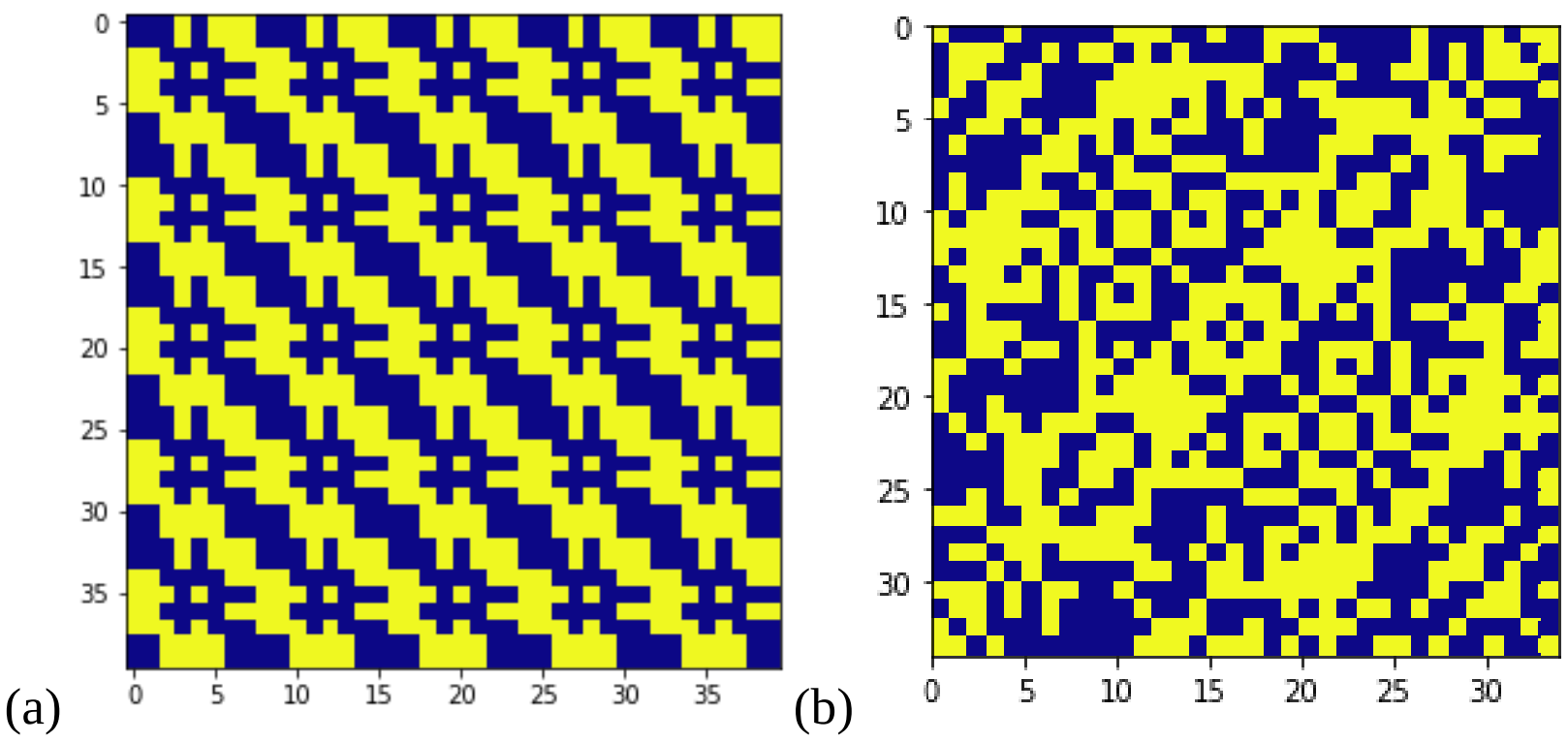

*Figure 1: (a) a distinct pattern of an 8x8 binary color-coded magic square repeated in a larger square grid and (b) a binary color-coded magic square of order 35.*

## Distinct Binary Magic Squares: A PCA analysis

PCA (Principal Component Analysis) [13] is a technique used to reduce the dimensionality of a dataset. In a Python code (supplementary data) we have developed (e.g. 4×4 magic squares), the PCA is applied to the flattened binary magic squares dataset to reduce the dimensionality from 16 to 2. The PCA finds the principal components of the dataset, which are linear combinations of the original features that capture the most variance in the data. The first principal component captures the most variance in the data, followed by the second principal component, and so on. By reducing the dimensionality of the data, PCA helps to simplify the dataset and make it easier to visualize.

The result (Figure 2) shows a scatter plot of the binary magic squares in 2-dimensional space after performing PCA on the flattened binary matrix of each square. Each point in the scatter plot represents a binary magic square, and the color of the point represent the pattern string of the square. The scatter plot can be used to interpret the relationships between the binary magic squares in the dataset. Points that are close together in the plot represent binary magic squares that are similar to each other in terms of their binary patterns. Similarly, points that are far apart from each other in the plot represent binary magic squares that are dissimilar to each other in terms of their binary patterns. The plot also shows that there are several clusters of points in the plot. Each cluster represents a group of binary magic squares that have similar binary patterns. By examining the plot and the pattern strings associated with each cluster, it is possible to identify the different types of patterns that occur in the binary magic squares.

Based on Figure 2, there are 24 distinct patterns among the 4×4 binary magic squares, 44 distinct patterns among the 5×5 associative magic squares, and 32 distinct patterns among the 8×8 Franklin magic squares, as identified through the use of PCA. Additionally, each plot displays a symmetrical distribution of points, which are categorized into 9 regions for 4×4 magic squares, 6 regions for 5×5 associative magic squares, and 9 regions for 8×8 Franklin magic squares. These regions represent different clusters of points in the scatter plot and are indicative of similarities in the binary patterns of the magic squares within each cluster.

Based on the PCA analysis, Table 2 presents the number of squares for each distinct pattern of different types of magic squares. It can be observed that the number of squares also exhibits symmetry.

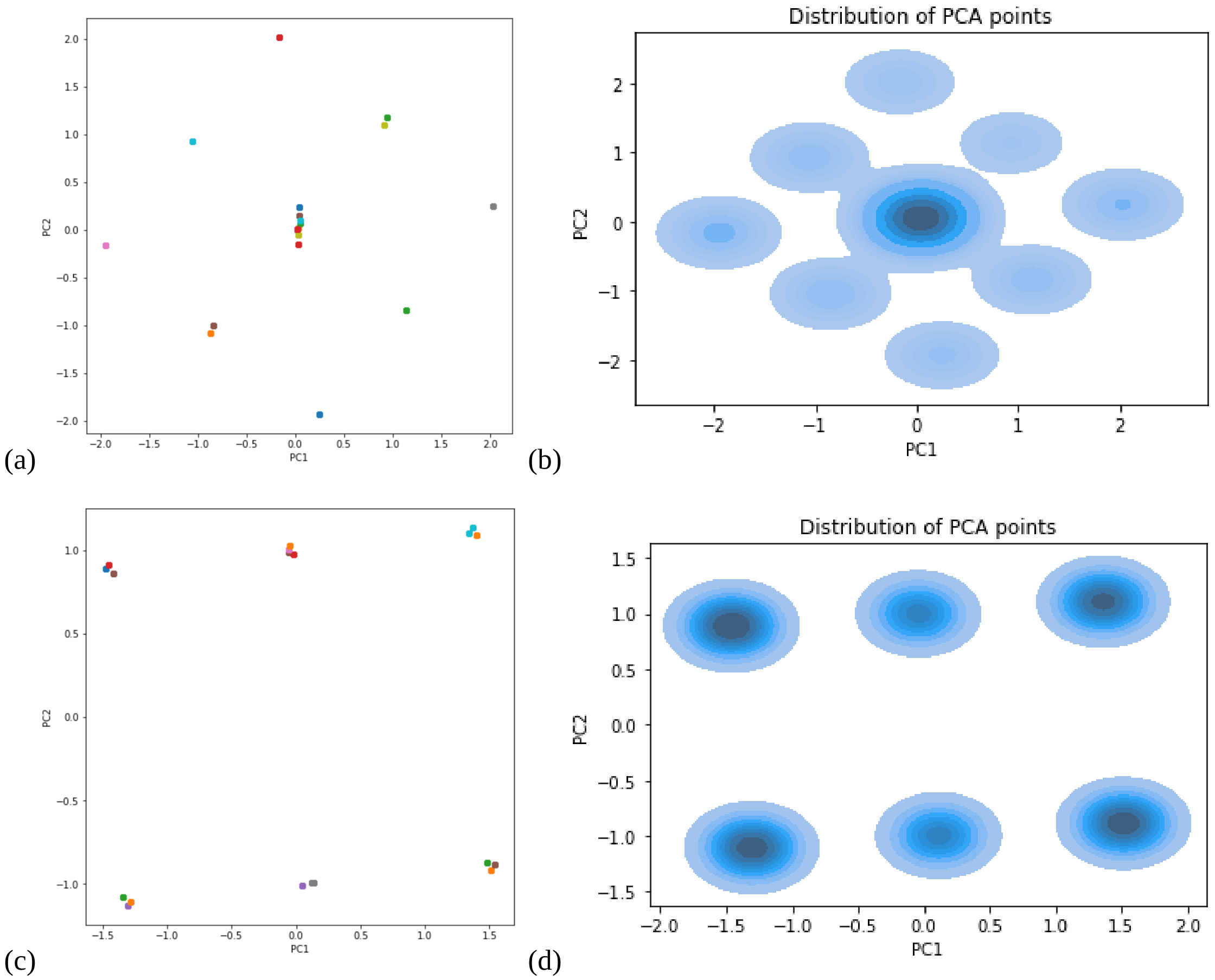

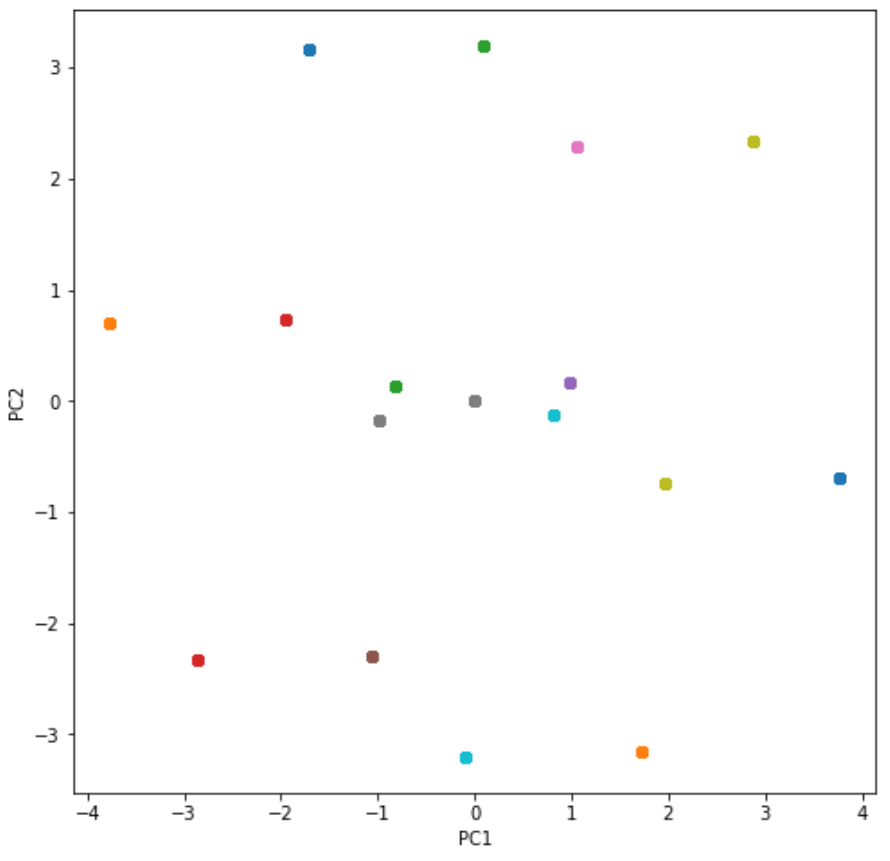

(e)

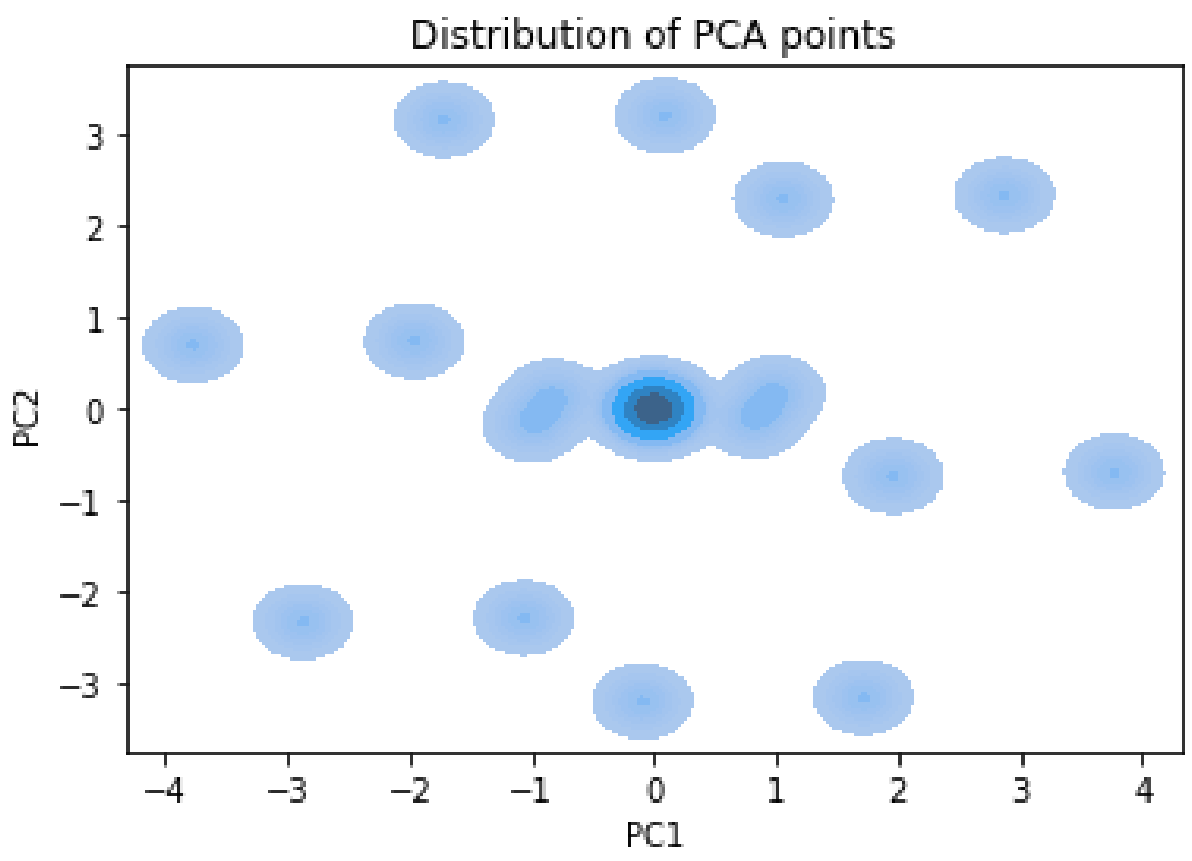

(f)

*Figure 2: PCA analysis of the binary 4×4 magic squares (a,b), 5×5 associative magic squares (c,d), and 8×8 Franklin magic squares (e,f).*

*Table 2: Number of binary squares with a distinct pattern based on PCA analysis.*

| **4×4 magic squares** | **5×5 associative magic squares** | **8×8 Franklin magic squares** |
|---|---|---|
| 0011001111001100: 24 | 0000101011111111101010000: 930 | 0011001111001100001100111100110000110011110011000011001111001100: 13824 |
| 0011010110101100: 24 | 0000111010111110101110000: 930 | 0011011011001001001101101100100100110110110010010011011011001001: 13824 |
| 0011101001011100: 15 | 0001010011111111100101000: 1011 | 0011100111000110001110011100011000111001110001100011100111000110: 13824 |
| 0011110000111100: 34 | 0001011001111111001101000: 1011 | 0011110011000011001111001100001100111100110000110011110011000011: 13824 |
| 0011110011000011: 35 | 0010001110111110111000100: 864 | 0101010101010101101010101010101001010101010101011010101010101010: 9216 |
| 0101001111001010: 27 | 0010010101111111010100100: 820 | 0101010101010101101010101010101001010101101010101010101001010101: 9216 |
| 0101010110101010: 62 | 0011101110001000111011100: 1016 | 0101010101010101101010101010101010101010010101010101010110101010: 9216 |
| 0101101001011010: 24 | 0011110101001001010111100: 892 | 0101010101010101101010101010101010101010101010100101010101010101: 9216 |
| 0101101010100101: 54 | 0100010011111111100100010: 1011 | 0101010110101010101010100101010101010101010101011010101010101010: 9216 |
| 0101110000111010: 24 | 0100011001111111001100010: 1011 | 0101010110101010101010100101010101010101101010101010101001010101: 9216 |
| 0110011010011001: 61 | 0101100001111111000011010: 1008 | 0101010110101010101010100101010110101010010101010101010110101010: 9216 |
| 0110100101101001: 42 | 0101110000111110000111010: 1008 | 0101010110101010101010100101010110101010101010100101010101010101: 9216 |
| 1001011010010110: 55 | 0101110011001001100111010: 962 | 0110001110011100011000111001110001100011100111000110001110011100: 13824 |
| 1001100101100110: 56 | 0101111001001001001111010: 962 | 0110011010011001011001101001100101100110100110010110011010011001: 13824 |
| 1010001111000101: 24 | 0110101110001000111010110: 1016 | 0110100110010110011010011001011001101001100101100110100110010110: 13824 |
| 1010010101011010: 61 | 0110110101001001010110110: 892 | 0110110010010011011011001001001101101100100100110110110010010011: 13824 |
| 1010010110100101: 20 | 0111000100111110010001110: 908 | 1001001101101100100100110110110010010011011011001001001101101100: 13824 |
| 1010101001010101: 55 | 0111000111001001110001110: 1262 | 1001011001101001100101100110100110010110011010011001011001101001: 13824 |
| 1010110000110101: 12 | 0111001101001001011001110: 1262 | 1001100101100110100110010110011010011001011001101001100101100110: 13824 |
| 1100001100111100: 60 | 0111010110001000110101110: 1262 | 1001110001100011100111000110001110011100011000111001110001100011: 13824 |
| 1100001111000011: 41 | 0111011100001000011101110: 1262 | 1010101001010101010101011010101001010101010101011010101010101010: 9216 |
| 1100010110100011: 26 | 1000001011111111101000001: 1266 | 1010101001010101010101011010101001010101101010101010101001010101: 9216 |
| 1100101001010011: 24 | 1000011010111110101100001: 1266 | 1010101001010101010101011010101010101010010101010101010110101010: 9216 |
| 1100110000110011: 20 | 1001100010111110100011001: 1193 | 1010101001010101010101011010101010101010101010100101010101010101: 9216 |
| | 1001101000111110001011001: 1193 | 1010101010101010010101010101010101010101010101011010101010101010: 9216 |
| | 1001101011001001101011001: 1090 | 1010101010101010010101010101010101010101101010101010101001010101: 9216 |
| | 1001111010001000101111001: 1090 | 1010101010101010010101010101010110101010010101010101010110101010: 9216 |
| | 1010100100111110010010101: 864 | 1010101010101010010101010101010110101010101010100101010101010101: 9216 |
| | 1010100111001001110010101: 1080 | 1100001100111100110000110011110011000011001111001100001100111100: 13824 |
| | 1010101101001001011010101: 1080 | 1100011000111001110001100011100111000110001110011100011000111001: 13824 |
| | 1010110110001000110110101: 1080 | 1100100100110110110010010011011011001001001101101100100100110110: 13824 |
| | 1010111100001000011110101: 1080 | 1100110000110011110011000011001111001100001100111100110000110011: 13824 |
| | 1011001110001000111001101: 1334 | |
| | 1011010101001001010101101: 1424 | |
| | 1100100010111110100010011: 1193 | |
| | 1100101000111110001010011: 1193 | |
| | 1100101011001001101010011: 1090 | |
| | 1100111010001000101110011: 1090 | |
| | 1101000001111111000001011: 1222 | |
| | 1101010000111110000101011: 1222 | |
| | 1101010011001001100101011: 1218 | |
| | 1101011001001001001101011: 1218 | |
| | 1110001110001000111000111: 1334 | |
| | 1110010101001001010100111: 1424 | |

**Distinct Binary Magic Squares: A LDA analysis**

Linear Discriminant Analysis (LDA) [14] is a statistical method used for classification and dimensionality reduction. It is based on finding the linear combination of features that maximally separates different classes in a dataset. LDA is a supervised learning algorithm, meaning that it requires labeled data in order to learn the optimal linear discriminant. The goal of LDA is to find a projection of the data that maximizes the ratio of between-class variance to within-class variance. In other words, it seeks a linear transformation that maximizes the distance between the means of the different classes, while minimizing the variance within each class. This is achieved by solving an eigenvalue problem involving the scatter matrices of the data.

In terms of binary magic squares, Linear Discriminant Analysis (LDA) is used to find the most important linear combination of the features (i.e., the binary values in each square) that can best separate the magic squares into different classes. In the python code (supplementary data) we have developed (for example, for 4×4 magic squares), LDA is used to find a linear transformation of the original 16-dimensional feature space to a 2-dimensional space that maximally separates the different patterns of the binary magic squares. The data points in the LDA analysis of binary magic squares provide information about how different magic squares are related to each other based on their binary patterns. Each data point (Figure 3) corresponds to a unique binary pattern, and the position of the data point in the LDA plot reflects how similar or dissimilar that pattern is to the other patterns in the dataset. By using LDA to analyze the binary patterns of magic squares, we can identify clusters of squares that are similar to each other in terms of their patterns.

Figure 3a depicts the classification of 4×4 binary magic squares into 7 pencil-like fingerprints based on LDA analysis, which reveals a distinct pattern. The colors in Figure 3a correspond to the LDA data points, which are based on the unique reflection/rotation binary patterns (Table 1). The distribution of LDA values is visualized in the form of a heatmap and a histogram in Figures 3b to 3d. Additionally, the LDA analysis results for 5×5 associative magic squares are shown in Figures 3e to 3h. The histograms of LDA values exhibit an approximate normal distribution. In Figures 3g and 3h, we compare the distributions of LDA1 and LDA2 with a normal distribution to assess their similarity.

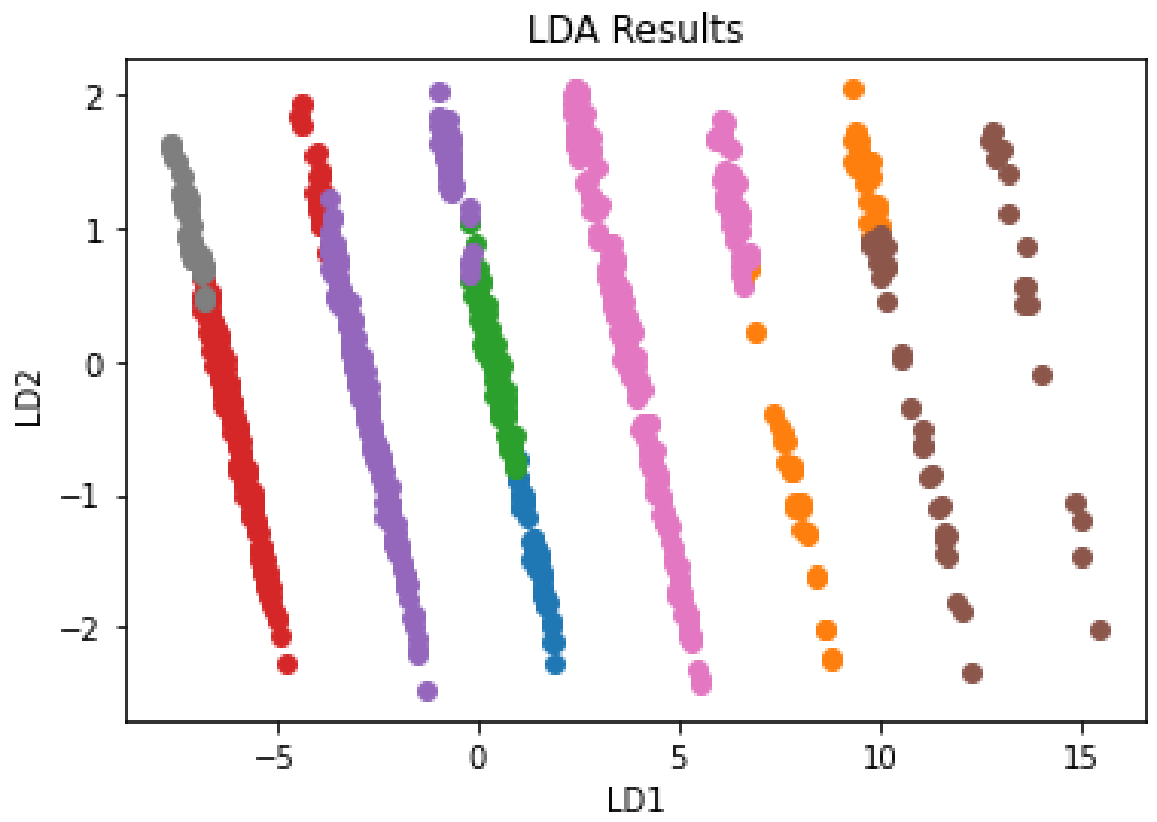

(a)

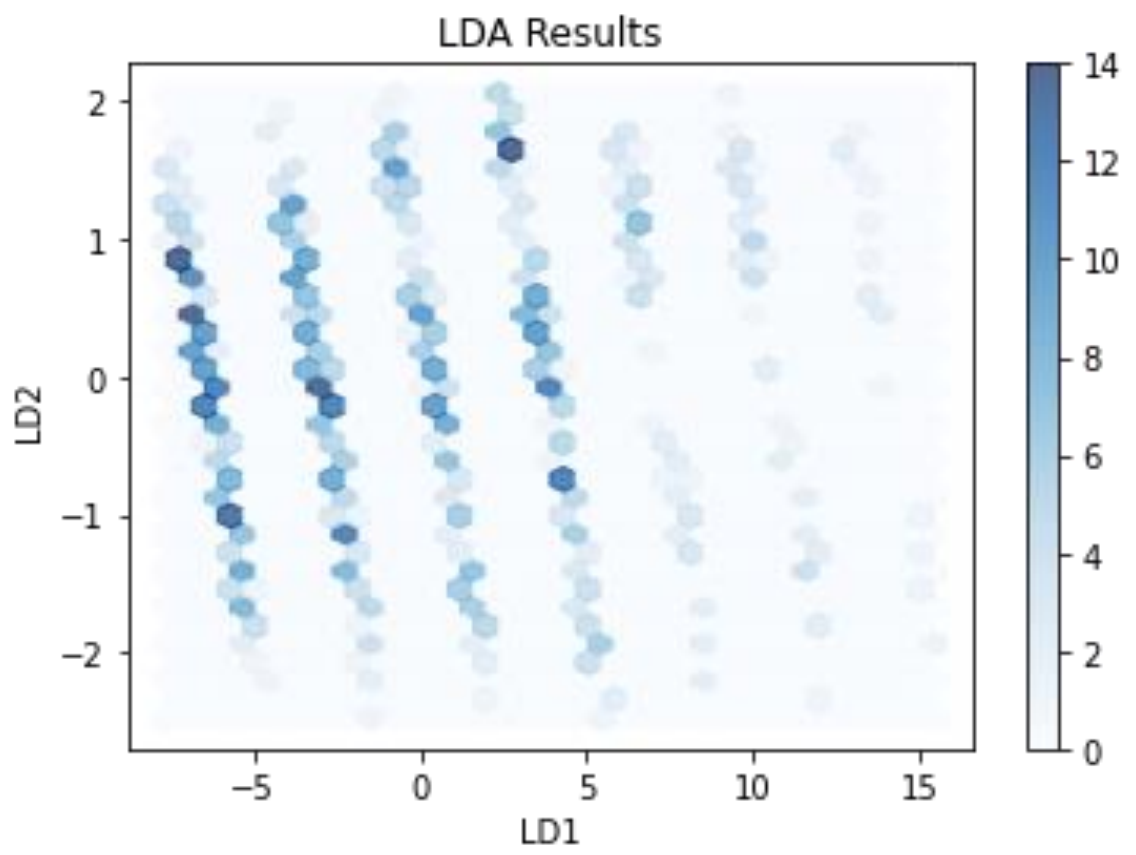

(b)

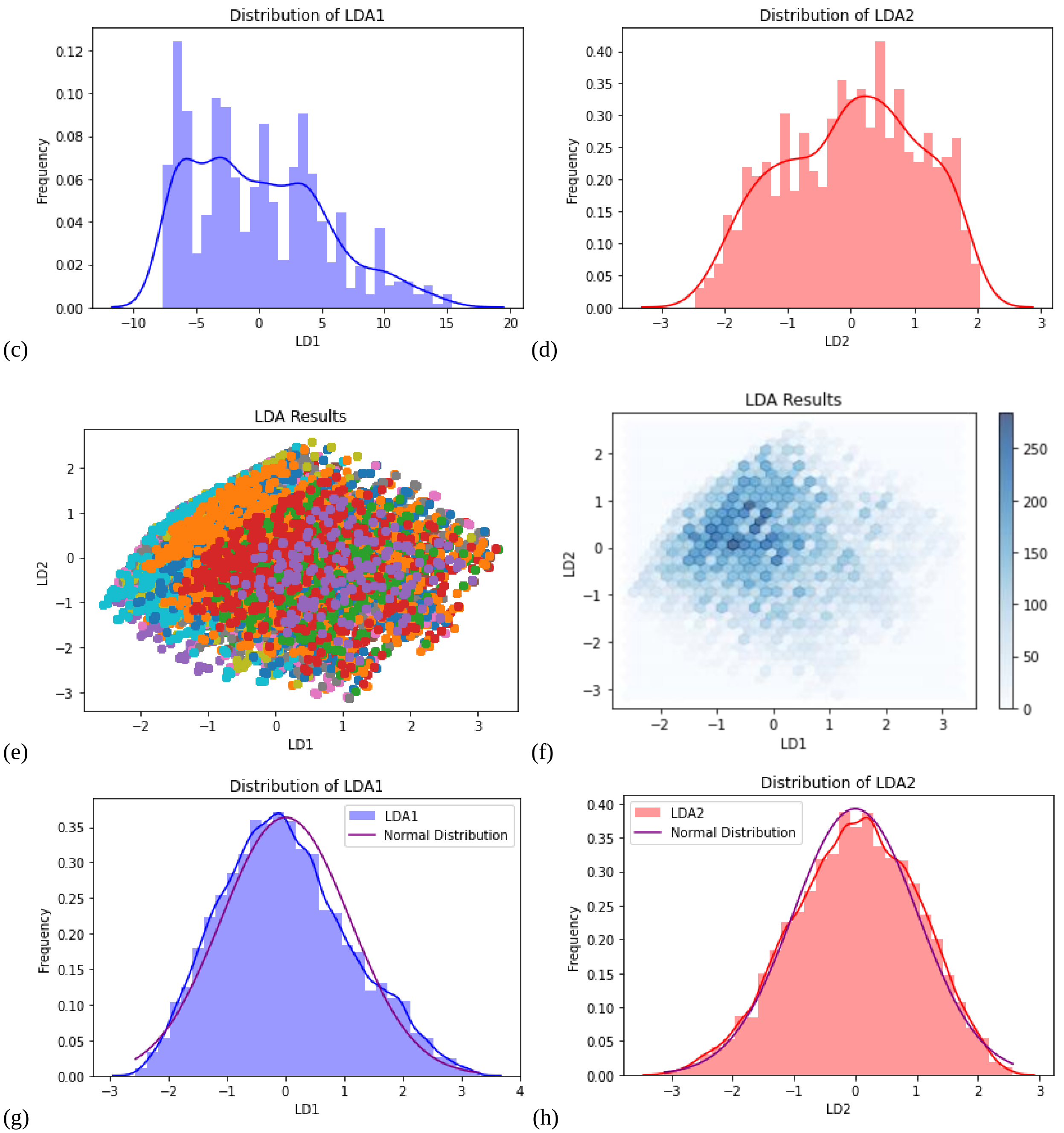

*Figure 3: Distinct patterns of binary 4×4 magic squares (a-d) and 5×5 associative magic square (e-h) based on the LDA analysis.*

## Conclusion

In this study, we have introduced the concept of binary color-coded magic squares, a novel approach that reveals previously unknown patterns and relationships within magic squares. We have explored the uniqueness of patterns within these squares using rotation and reflection, PCA, and LDA. Our investigation has revealed striking patterns that were previously unknown in traditional magic squares. Based on the analysis conducted in this study, the following results can be concluded:

- By analyzing binary color-coded magic squares through rotation and reflection, we identified the distinct color-coded patterns and their frequency in different types of magic squares. We found that only a small number of unique binary matrices are needed to construct a large number of magic squares.
- Principal component analysis (PCA) allowed us to identify the most significant patterns within the dataset of binary color-coded magic squares. The scatter plot resulting from PCA showed that similar binary magic squares are grouped together, suggesting underlying patterns and relationships.
- Linear discriminant analysis (LDA) was able to classify and distinguish between different types of magic squares based on their attributes, suggesting that the binary color-coding element provides additional information that can be used to differentiate between different types of magic squares.

## Acknowledgment

I would like to express my appreciation to Professor Chérif F. Matta (Mount Saint Vincent University) and Professor Thanh-Tung Nguyen-Dang (Université Laval) for their support and guidance throughout this study. Also, I would like to extend my thanks to Craig Knecht and Walter Trump (Gymnasium Stein) for providing me with useful information that enabled me to access data from various websites. Finally, I am grateful to the anonymous reviewers for their comments, which helped in enhancing the quality of this work.

## Conflict of Interest

The author has no conflicts of interest to disclose.

## Supplementary Data

The complete Python code used for generating the results in Tables 1 and 2, and Figures 1 to 3, is provided, along with step-by-step explanations of the methodology. Readers can access this data online for transparency and reproducibility.